\documentclass[a4paper,11pt]{article}
\usepackage{amssymb}
\usepackage{latexsym}
\newtheorem{theorem}{Theorem}[section]
\newtheorem{lemma}[theorem]{Lemma}

\newtheorem{definition}[theorem]{Definition}
\newtheorem{proposition}[theorem]{Proposition}

\newtheorem{remark}[theorem]{Remark}

\newtheorem{conj}[theorem]{Conjecture}

\newcommand{\comment}[1]{}  

\newcommand{\res}{\restriction}

\newcommand{\CZF}{{\mathbf{CZF}}}

\newcommand{\ZF}{{\mathbf{ZF}}}

\newcommand{\IKP}{{\mathbf{IKP}}}

\newcommand{\AC}{{\mathrm{AC}}}

\newcommand{\prf}{{\bf Proof: }}

\newcommand{\beqs}{\begin{eqnarray*}}
\newcommand{\eeqs}{\end{eqnarray*}}

\newcommand{\KP}{{\mathbf{KP}}}

\newcommand{\beq}{\begin{eqnarray}}
\newcommand{\eeq}{\end{eqnarray}}

\newcommand{\IZF}{{\mathbf{IZF}}}

\newcommand{\dom}{{\mathrm{dom}}}

\newcommand{\PAs}[1]
{#1^+}

\newcommand{\bes}{\begin{eqnarray*}}
\newcommand{\ees}{\end{eqnarray*}}

\def\provx#1#2#3#4{
\setbox1=\hbox{\kern1.5pt$\scriptstyle#3$}
\def\zeichen{#2}
\ifx\zeichen\empty\setbox0=\hbox to .75em{}\else\setbox0=\hbox
{\kern1.5pt$\scriptstyle#2$}\fi
\dimen1=\dp0 \ifdim \dimen1=0pt
\advance \dimen1 by 1.5ex \else \advance \dimen1 by 1.2ex
\fi\dimen3=2ex\dimen4=.5ex\ifdim \wd0<\wd1 \dimen2=\wd1 \else \dimen2=\wd0
\fi\hbox{$#1\hskip 5pt minus5pt\vrule height\dimen3
depth\dimen4\raise\dimen1\copy0\hskip-1\wd0 \lower\ht1
\copy1\hskip-1\wd1\vrule width\dimen2 height.7ex depth-.6ex\hskip3pt
minus1.5pt#4\hskip2pt plus2pt minus2pt$}}

\def\prov#1#2#3{
\setbox1=\hbox{\kern1.5pt$\scriptstyle#2$}
\def\zeichen{#1}
\ifx\zeichen\empty\setbox0=\hbox to .75em{}\else\setbox0=\hbox
{\kern1.5pt$\scriptstyle#1$}\fi
\dimen1=\dp0
\ifdim \dimen1=0pt
\advance \dimen1 by 1.5ex \else \advance \dimen1 by 1.2ex
\fi\dimen3=2ex\dimen4=.5ex\ifdim \wd0<\wd1 \dimen2=\wd1 \else \dimen2=\wd0
\fi\hbox{\hskip0pt plus 4pt
$\vrule height\dimen3
depth\dimen4\raise\dimen1\copy0\hskip-1\wd0
\lower\ht1\copy1\hskip-1\wd1\vrule width\dimen2 height.7ex depth-.6ex
\hskip3pt minus1.5pt#3\hskip2pt plus2pt minus2pt$}}

\def\prv#1#2{
\setbox1=\hbox{\kern1.5pt$\scriptstyle#2$}
\ifx\zeichen\empty\setbox0=\hbox to .75em{}\else\setbox0=\hbox
{\kern1.5pt$\scriptstyle#1$}\fi
\dimen1=\dp0 \ifdim \dimen1=0pt
\advance \dimen1 by 1.5ex \else \advance \dimen1 by 1.2ex
\fi\dimen3=2ex\dimen4=.5ex\ifdim \wd0<\wd1 \dimen2=\wd1 \else \dimen2=\wd0
\fi\hbox{\hskip.5em$\vrule height\dimen3
depth\dimen4\raise\dimen1\copy0\hskip-1\wd0
\lower\ht1\copy1\hskip-1\wd1\vrule width\dimen2 height.7ex depth-.6ex
\hskip3pt minus1.5pt$}}

\mathchardef\str='1066
\def\negprov#1#2#3{
\setbox1=\hbox{\kern1.5pt$\scriptstyle#2$}
\setbox4=\hbox{$\str$}
\def\zeichen{#1}
\ifx\zeichen\empty\setbox0=\hbox to 1em{}\else\setbox0=\hbox
{\kern1.5pt$\scriptstyle#1$}\fi
\dimen1=\dp0
\ifdim \dimen1=0pt
\advance \dimen1 by 1.5ex \else \advance \dimen1 by 1.2ex
\fi\dimen3=2ex\dimen4=.5ex\ifdim \wd0<\wd1 \dimen2=\wd1 \else \dimen2=\wd0
\fi\hbox{\hskip.5em$\kern-1.9pt\raise1pt\copy4\kern-\wd4\kern1.9pt\vrule height\dimen3
depth\dimen4\raise\dimen1\copy0\hskip-1\wd0
\lower\ht1\copy1\hskip-1\wd1\vrule width\dimen2 height.7ex depth-.6ex
\hskip3pt minus1.5pt#3\hskip2pt plus2pt minus2pt$}}

\def\goed#1{\setbox5=\hbox{$#1$}\dimen1=.25em \dimen2=\dimen1 \advance \dimen2
by -1pt\hbox{\raise.65\ht5 \hbox{\vrule height.5\ht5 depth0pt width.4pt\vrule
height.5\ht5 width\dimen1 depth-.48\ht5}\kern-\dimen2\copy5\kern-\dimen2
\raise.65\ht5 \hbox{\vrule height .5\ht5 width\dimen1 depth-.48\ht5\vrule
height.5\ht5 depth 0pt width.4pt}\hskip4pt plus2pt minus2pt}}

\def\mod#1#2{
\def\zeichen{#1}
\hbox{\hskip 2pt plus3pt minus 2pt\vrule width.5pt height2ex depth.5ex
\vbox{\ifx\zeichen\empty\hbox to .75em{}\else
\hbox{\kern1.5pt $\scriptstyle#1$}\fi
\kern2pt
\hrule
\kern1.7pt
\hrule\kern1.7pt}
\hskip3pt minus 2pt$#2$}\hskip2pt
plus3pt minus2pt}

\def\notmod#1#2{\hbox{\hskip 2pt plus 3pt minus 3pt\vrule width.5pt
height2ex depth.5ex
\vbox{\hbox{\kern1.5pt $\scriptstyle#1$}\kern3pt
\setbox0=\hbox{\kern2pt$\scriptstyle/$}
\hrule
\kern-1.7pt
\copy0
\kern-\ht0
\kern 1.7pt
\hrule\kern1.7pt}n
\hskip3pt minus 2pt$#2$}\hskip2pt
plus3pt minus2pt}
\def\sq{\hbox{\rlap{$\sqcap$}$\sqcup$}}
\def\qed{\ifmmode\sq\else{\unskip\nobreak\hfil\penalty50\hskip1em\null
\nobreak\hfil\sq\parfillskip=0pt\finalhyphendemerits=0\endgraf}\fi\medskip}

\def\lleq{\hbox{\hskip3pt minus3pt\kern1pt\lower4pt
\vbox{\hbox{$\scriptstyle\ll$}
\kern-7pt\hbox{\kern1pt$\scriptstyle=$}}\hskip3pt minus 3pt}}

\mathchardef\res='1152
\mathchardef\qin='1062
\mathchardef\qprec='1036
\mathchardef\qless='474
\mathchardef\dpkt='72

\newcommand{\CH}{\mathrm{CH}}
\newcommand{\Vdasha}{\Vdash_{\!\mbox{\tiny $A$}}}
\newcommand{\Vdashas}{\Vdash_{\!\mbox{\tiny $A'$}}}

\newcommand{\Vdashce}{\Vdash_{\!\mbox{\tiny $C\cup E$}}}
\newcommand{\hocha}{^{\mbox{\tiny $L[A]$}}}
\newcommand{\hochc}{^{\mbox{\tiny $L[C]$}}}
\newcommand{\hochas}{^{\mbox{\tiny $L[A']$}}}

\newcommand{\hochal}{^{\mbox{\tiny $L_{\lambda}[A]$}}}

\begin{document}

\title{Indefiniteness in semi-intuitionistic set theories: On a conjecture of Feferman}

\author{Michael Rathjen\\
Department of Pure Mathematics\\
University of Leeds, Leeds LS2 9JT, England\\  E-mail:~{\sf
rathjen@maths.leeds.ac.uk} }
\maketitle
\begin{abstract}
The paper proves a conjecture of Solomon Feferman concerning the indefiniteness of the continuum
hypothesis relative to a semi-intuitionistic set theory.
\\[2ex]

\end{abstract}

\section{Introduction}
Frege in \cite[$\S$ 68]{frege} wrote: {\em Ich setze voraus,
dass man wisse, was der Umfang eines Begriffes sei.}\footnote{Translation: {\em I assume that it is known what the extension of a concept is.}} Dummett's diagnosis of the failure of Frege's logicist project in the final chapter of \cite{dummett} focusses on the adoption of classical quantification over domains comprised of objects
 falling under an indefinitely extensible concept.
He repudiates the classical view as illegitimate and puts forward reasons  in favor of an intuitionistic interpretation of quantification.
Solomon Feferman, in recent years, has argued that the Continuum Hypothesis ($\CH$) might  not be a definite
mathematical problem (see \cite{F2,F3,F4}\footnote{Incidentally, the paper \cite{F2} was written for Peter Koellner's
 {\em Exploring the frontiers of incompleteness} (EFI) Project, Harvard 2011-2012.}).
\begin{quote} My reason for that is that the concept of arbitrary set essential to its formulation is vague or underdetermined and there is no way to sharpen it without violating what it is supposed to be about.
In addition, there is considerable circumstantial evidence to support the view that $\CH$ is not definite.
(\cite[p.1]{F2}). \end{quote}
In  particular the power set, ${\mathcal P}(A)$, of a given set $A$ may be considered to be an indefinite collection whose members are subsets of $A$, but whose exact extent is indeterminate (open-ended).
In \cite{F2}, Feferman proposed a logical framework for what's definite and for what's not.
\begin{quote} One
way of saying of a statement $\varphi$ that it is definite is that it is true or false; on a deflationary
account of truth that's the same as saying that the Law of Excluded Middle (LEM) holds
of $\varphi$ , i.e. one has $\varphi\,\vee\,\neg\varphi$ . Since LEM is rejected in intuitionistic logic as a basic principle,
that suggests the slogan, ``What's definite is the domain of classical logic, what's not is
that of intuitionistic logic.'' [...] And in the case of set theory, where every set is
conceived to be a definite totality, we would have classical logic for bounded
quantification while intuitionistic logic is to be used for unbounded quantification.
(\cite[p. 23]{F2})\end{quote}
At the end of \cite{F2} he made that idea more precise by suggesting semi-intuitionistic set theories as frameworks for formulating questions of definiteness and studying the definiteness of specific set-theoretic statements.
In relation to $\CH$, he conjectured that this statement is not definite in the specific case of a semi-intuitionistic set theory $\mathbf T$, in the sense that $\mathbf T$ does not prove $\CH\,\vee\,\neg \CH$. The set-theoretical point of view expressed by $\mathbf T$ accepts the definiteness of the continuum
in its guise as the arithmetical/geometric structure of the real line, but does not allow the powerset operation to be applied to arbitrary sets.

The objective of this paper is to  prove Feferman's conjecture. In this sense it is a technical paper.
It lays out new evidence for the reader to consider. However,
as far as the ongoing discussions of the foundational status of $\CH$ are concerned, readers will have to form their own conclusions.

A chief technique applied in this article is realizability over relativized constructible hierarchies
 combined with forcing. More widely the impression is that $\CH$ is not an isolated case in that other statements could be proved to be indefinite relative to semi-intuitionistic set theories in this way.
At any rate, it appears that the paper adds a hitherto unexplored tool to the weaponry earmarked for engineering specific realizability models
and proving independence results.

An outline of the paper of the paper reads as follows: Section 2 introduces formal systems of semi-intuitionistic set theory and in particular the theory $\mathbf T$. Section 3 is devoted to the relativized constructible hierarchy $L[A]$ and
its properties. In section 4, $L[A]$ features as a domain of computation which gets utilized in
section 5 as a realizability universe for $\mathbf T$. By carefully designing sets of ordinals $C$ and $E$
and employing results from forcing, realizability of $\mathbf T$ over $L[C]$ and $L[C\cup E]$ yields
conflicting information that leads to a contradiction, and thus provides a proof of
the desired conjecture.

\section{Semi-intuitionistic set theory}
 The
study of subsystems of $\ZF$ formulated in intuitionistic logic with Bounded    Separation
 was apparently initiated by Pozsgay \cite{pozsgay1,pozsgay2} and then pursued more
systematically by Tharp \cite{tharp}, Friedman \cite{frieda} and Wolf \cite{wolf}.
These systems are actually semi-intuitionistic as they contain the law of excluded middle for bounded formulae.

Classical Kripke-Platek set theory, $\KP$, is an important theory that accommodates a great deal of set theory. Its transitive models, called
admissible sets, have been a major source of interaction between
model theory, recursion theory and set theory (cf. \cite{barwise}).
 $\KP$ arises from $\ZF$ by completely omitting the power set axiom and
restricting separation and collection to bounded formulae.
Here we are interested in its intuitionistic cousin.
\begin{definition}\label{def1}{\em {\em Intuitionistic Kripke-Platek set theory}, $\IKP$,
 is formulated in the
usual language of set theory
containing $\in$ as the only non-logical symbol besides $=$.
Formulae are built from prime formulae $a\in b$ and $a=b$ by use
of propositional connectives and quantifiers $\forall x,\exists
x$. Quantifiers of the forms $\forall x\in a$, $\exists x\in a$
are called {\em bounded}. {\em Bounded} or  {\em
$\Delta_0$-formulae} are the formulae wherein all quantifiers are
bounded.
$\IKP$ is based on
 intuitionistic logic and
has the following non-logical axioms:
{\em Extensionality, Pair, Union, Infinity} (in the specific version that there is a smallest set containing the empty set 0 and closed under the successor operation, $x'=x\cup\{x\}$), {\em Bounded Separation}
 $$\exists x\,\forall u\left[u\in x\leftrightarrow(u\in
 a\,\wedge\,\varphi(u))\right]$$
 for all bounded formulae $\varphi(u)$, {\em Bounded Collection}
 $$\forall x\in a\,\exists y\,\psi(x,y)\,\to\,\exists z\,\forall x\in
 a\,\exists y\in z\,\psi(x,y)$$
 for all bounded formulae $\psi(x,y)$, and
 {\em Set Induction}
 $$\forall x\,\left[(\forall y\in x\,\theta(y))\to
 \theta(x)\right]\to\,\forall x\,\theta(x)$$
 for all formulae $\theta(x)$. \\[1ex]
  Feferman in \cite{F1} proceeded to add several further schemata to the axioms of $\IKP$. The most basic principle that he added follows from the idea that in semi-constructive set theory each set is considered
  to be a definite totality. As a consequence of $\Delta_0$ separation
    one obtains
  a restricted Law of Excluded Middle:
   $$ ({\Delta_0}\mathrm{-LEM})\;\;\;\varphi\vee\neg\varphi, \mbox{
    for all $\Delta_0$-formulae $\varphi$.}$$
   {\em Markov's principle} in the form
   $$(\mathrm{MP})\;\;\;\neg\neg \exists x\varphi \to \exists x\varphi,\mbox{ for all $\Delta_0$ formulae $\varphi$}$$
   is another principle that is frequently added in this context.

   Some further principles that are considered in \cite{F1} are
   $(\mathrm{BOS})$ and $\mathrm{AC}_{\mathrm{Set}}$.
    $$(\mathrm{BOS})\;\;\;
\forall x\in a\,[\varphi(x)\vee\neg\varphi(x)]\to
 [\forall x\in a\,\varphi(x) \;\vee\; \exists x\in a\,\neg\varphi(x)]$$
 for all formulae $\varphi(x)$.
 $$(\mathrm{AC}_{\mathrm{Set}})\;\;\;
 \forall x\in a\,\exists y\,\psi(x,y) \to \exists f\,[\mathrm{Fun}(f)\,\wedge\,\mathrm{dom}(f)=a\,\wedge\,\forall
 x\in a\,\varphi(x,f(x))]$$
 for all formulae $\psi(x,y)$, where $\mathrm{Fun}(f)$ expresses in the usual set-theoretic form that $f$ is a function, and $\mathrm{dom}(f)=a$ expresses that
 the domain of $f$ is the set $a$.
} \end{definition}
Feferman \cite[Theorem 6]{F1} shows that $\mathbf{SCS}:=\IKP+({\Delta_0}\mathrm{-LEM})+(\mathrm{MP})+(\mathrm{BOS})+
(\mathrm{AC}_{\mathrm{Set}})$ has the same proof-theoretic strength as $\KP$
(and therefore the same as $\IKP$). His proof uses a functional interpretation.
The same result can be obtained via a realizability interpretation using codes for
$\Sigma_1$ partial recursive set functions as realizers along the lines
of Tharp's 1971 arcticle \cite{tharp}.

\begin{remark} {\em \begin{itemize}
\item[(i)] $\mathbf{SCS}$ proves the full replacement schema of $\ZF$. Moreover,
$\mathbf{SCS}$ proves strong collection, i.e. all formulae
$$\forall x\in a\,\exists y\,\varphi(x,y) \to \exists z\,[\forall x\in a \,\exists y\in
z\,\varphi(x,y)\,\wedge\,\forall y\in z\,\exists x\in a\,\varphi(x,y)]$$
where $\varphi(x,y)$ is an arbitrary formula.

 Strong collection is an axiom schema of
Constructive Zermelo-Fraenkel set theory, $\CZF$ (cf. \cite{maralt,mar}) and also
of Tharp's set theory \cite{tharp}.
\item[(ii)] $\mathbf{SCS}$ is a subtheory of Tharp's semi-intuitionistic set theory $\IZF$ \cite{tharp}, for if
$\forall x\in a\,\exists y\,\varphi(x,y)$ holds, then there is a set
$d$ such that $\forall x\in a\,\exists z\in d\,\exists y[z=\langle x,y\rangle\,\wedge\,\varphi(x,y)]$ and $\forall z\in d\,\exists x\in a
\,\exists y[z=\langle x,y\rangle\,\wedge\,\varphi(x,y)]$
 using (strong) collection,
and by axiom 6 of $\IZF$, $d$ is the surjective image of an ordinal, i.e., there is
an ordinal $\alpha$ and a function $g$ with domain $\alpha$ and range $d$.
Thus $d$ is a set of ordered pairs.
Now define a function $f$ with domain $a$ by letting $f(x)$ be the second projection of $g(\xi)$ where
$\xi$ is the least ordinal $<\alpha$ such that the first projection of $g(\xi)$ equals
$x$.\end{itemize}
}
\end{remark}

As it turns out, some of the axioms of $\mathbf{SCS}$ are redundant.
\begin{proposition} $\IKP+(\mathrm{AC}_{\mathrm{Set}})$ proves $({\Delta_0}\mathrm{-LEM})$
and $(\mathrm{BOS})$. \end{proposition}
\prf First we prove $({\Delta_0}\mathrm{-LEM})$, using Diaconescu's old constructions
\cite{Diaconescu}.
Let $0$ be the empty set, $1:=\{0\}$ and $A=\{0,1\}$.
Note that (intuitionistically) $\forall x,y\in A\,[x=y\vee x\ne y]$
(where $x\ne y$ abbreviates $\neg x=y$) since $0\ne 1$ as $0\in 1$ and $0\notin 0$.
Suppose $\varphi$ is $\Delta_0$. Define $a:=\{n\in A\mid n=0\,\vee\,[n=1\wedge \varphi]\}$ and $b:=\{n\in A\mid n=1\,\vee\,[n=0\wedge \varphi]\}$. $a$ and $b$ are sets by $\Delta_0$ separation. Obviously we have $$\forall z\in \{a,b\}\,\exists
k\in A\;k\in z$$ since $0\in a$ and $1\in b$. So
we may apply $(\mathrm{AC}_{\mathrm{Set}})$ to obtain a function $f$ with domain $\{a,b\}$ such that $f(a),f(b)\in A$. We thus have
$f(a)=f(b)$ or $\neg(f(a)=f(b))$. In the first case, we can infer that $\phi$.
In the second case, we have $a\ne b$. As $\varphi$ implies $a=b$, we get $\neg \varphi$.

To show $(\mathrm{BOS})$ assume $$\forall x\in a \,[\psi(x)\,\vee\,\neg\psi(x)],$$
where $\psi(x)$ is an arbitrary formula. Thus,
$$\forall x\in a\,\exists y\,[(\psi(x)\wedge y=0)\,\vee\,(\neg\psi(x)\wedge y=1)].$$
With the help of $(\mathrm{AC}_{\mathrm{Set}})$ there is a function $f$ with domain
$a$ such that
\begin{eqnarray}\label{prop1}&&\forall x\in a\,[(\psi(x)\wedge f(x)=0)\,\vee\,(\neg\psi(x)\wedge f(x)=1)]\end{eqnarray}
and hence $\forall x\in a\,[ f(x)=0)\,\vee\, f(x)=1].$
Using $({\Delta_0}\mathrm{-LEM})$ we have $\exists x\in a\, f(x)=1$ or
$\forall x\in a\,f(x)=0$. In the former case we deduce $\exists x\in a\,\neg\psi(x)$ from
(\ref{prop1}), whereas in the latter case we infer that $\forall x\in a\,\psi(x)$.
\qed

\begin{definition}{\em Let $\mathbf{T}$ be the theory $$\mathbf{SCS}+\mbox{`$\mathbb R$ is
a set'}$$ where $\mathbf{SCS}$ is from Definition \ref{def1} and
`$\mathbb R$ is
a set' asserts that the reals, $\mathbb R$, form a set. Since $\mathbf{SCS}$ has classical logic for $\Delta_0$ formulae it is not necessary to
 pay much attention to the question of how the reals are actually formalized as is so often the case in intuitionistic contexts. Thus,  any of the following equivalent statements could be used to formalize the existence of $\mathbb R$ as a set:
\begin{itemize}
\item
The collection of all functions from $\mathbb N$ to $\mathbb N$, ${\mathbb N}^{\mathbb N}$, is a set.
\item The collection of all subsets of $\mathbb N$ is a set.
 \end{itemize}
 }\end{definition}

 \begin{remark}{\em The proof-theoretic strength of $\mathbf T$ resides strictly between
 full classical second order arithmetic and Zermelo set theory. In particular all
 theorems of classical second order arithmetic are theorems of $\mathbf T$ and
 all theorems of $\mathbf T$ are theorems of Zermelo set theory plus the axiom of choice.}
 \end{remark}

 The continuum hypothesis, $\CH$, is the statement that every infinite set of reals is
 either in one-one correspondence with $\mathbb N$ or with $\mathbb R$.
 More formally, this can be expressed as follows:
 $$\forall x\subseteq \mathbb R\,[x\ne\emptyset\to
 (\exists f\, f:\omega\twoheadrightarrow x\,\vee\,\exists f\,f:x\twoheadrightarrow
 {\mathbb R})]$$
 where $f:y\twoheadrightarrow z$ signifies that $f$ is a surjective function
 with domain $y$ and co-domain $z$.

 \begin{conj}[Feferman] $\mathbf T$ does not prove $\CH\,\vee\, \neg \CH$.
 \end{conj}

When one ponders how to prove the conjecture one of the first ideas that comes to mind is that intuitionistic set theories $S$ very often have the disjunction property, i.e.,
if $S\vdash \psi\vee \theta$ then $S\vdash \psi$ or $S\vdash \theta$ (cf. \cite{tklr,tklrac}). If this property held for $\mathbf T$ it would certainly settle
the conjecture in the affirmative. However, $\mathbf T$ being semi-intuitionistic,
the disjunction property does not hold for it.
The technique of realizability certainly springs to mind when tackling such problems and consequently one would like to show that there is a realizability interpretation of $\mathbf T$ that has no realizer for
$\CH\,\vee\,\neg \CH$. There are several essentially different forms of realizability for set theories to choose from (cf. \cite{beeson,RMC,frieda,mc86,myhill73,ober,R04,tklr,tklrac,schwach,tharp}). Moreover, what should the realizers be and how should the realizability universe be defined?

\section{The relativized constructible hierarchy}
Later we shall look at realizability interpretations in the relativized constructible hierarchy.
The latter comes in two versions: For a set $A$ we have $L(A)$ and $L[A]$.
$L(A)$ is the smallest inner model that contains $A$.  In $L(A)$, the transitive closure of $A$ is added at level 0 and for higher levels the definition is the same as for $L$, whereas in $L[A]$
$A$ acts as an additional predicate for defining sets. The two hierarchies can be quite different.
 E.g., in general $L(A)$ is not a model of the axiom of choice, $\mathrm{AC}$, whereas
$L[A]$ is always a model of $\mathrm{AC}$. Another difference is that $L[{\mathbb R}]=L$
whereas  $L\ne L(\mathbb R)$ when ${\mathbb R}\notin L$.\footnote{Note that in the buildup of $L[{\mathbb R}]$, ${\mathbb R}$ is just used as a predicate. By identifying ${\mathbb R}$ with the set of all functions from $\mathbb N$ to $\mathbb N$, this is merely the predicate of being such a function, which is $\Delta_0$ in $\mathbb N$, hence absolute. Thus nothing outside of $L$ can be generated in this way.}
 Only $L[A]$ is interesting for the  purposes of this paper.

\begin{definition}{\em Let ${\mathcal L}_{\in}$ be the language of set theory
and ${\mathcal L}_{\in}(P)$ be its augmentation by a unary predicate symbol $P$.
Let $A$ be a set. Any set $X$ gives rise to a structure
$\langle X,\in,A\cap X\rangle$ for ${\mathcal L}_{\in}(P)$ with domain $X$ where the elementhood symbol is interpreted by the elementhood relation
restricted to $X\times X$ and $P$ is interpreted as $A\cap X$.
Thereby $A$ acts as a unary predicate on $X$.
 A subset $Y$ of $X$ is said to be definable in $\langle X,\in,A\cap X\rangle$
 if there is a formula $\varphi(x,y_1,\ldots,y_r)$ of ${\mathcal L}_{\in}(P)$
 with all free variables exhibited and $b_1,\ldots,b_r\in X$ such that for all $a\in X$,
 $$a\in Y\;\;\;\mbox{ iff }\;\;\;\langle X,\in,A\cap X\rangle\models\varphi(a,b_1,\ldots,b_r),$$
 where of course $\langle X,\in,A\cap X\rangle\models\varphi(a,b_1,\ldots,b_r)$ signifies that $\varphi$ holds in the structure under the variable assignment
 $x\mapsto a$ and $y_i\mapsto b_i$.

The sets $L_{\alpha}[A]$ are defined by recursion on $\alpha$ as follows:

 \begin{itemize}
 \item[(i)] ${\mathrm{Def}}^A(X)\,:=\,\{Y\subseteq X\mid\mbox{$Y$ definable in $\langle X,\in,A\cap X\rangle$}\}$.
     \item[(ii)] $L_0[A]=\emptyset$.
     \item[(iii)] $L_{\alpha+1}[A]={\mathrm{Def}}^A(L_{\alpha}[A])$.
     \item[(iv)] $L_{\lambda}=\bigcup_{\xi<\lambda}L_{\xi}[A]$ for limits $\lambda$.
     \item[(v)] $L[A]=\bigcup_{\alpha}L_{\alpha}[A]$.
 \end{itemize}
 }\end{definition}
 The next proposition lists some important properties of $L_{\alpha}[A]$.
 Bounded quantifiers are of the form $\forall x\in y$ and $\exists x\in y$.
 A bounded or $\Delta_0$ formula of ${\mathcal L}_{\in}(P)$ is a formula in which all quantifiers appear  bounded. A formula of ${\mathcal L}_{\in}(P)$ of the form $\exists z\varphi(z)$ ($\forall z\varphi(z)$) with $\varphi$ bounded is said to be $\Sigma_1$ ($\Pi_1$).
 Let $\alpha>0$. A relation on $L_{\alpha}[A]$ is said to be $\Sigma_1^{L_{\alpha}[A]}$
 ($\Pi_1^{L_{\alpha}[A]}$) if it is definable (with parameters) on the structure $\langle L_{\alpha}[A],\in,A\cap L_{\alpha}[A]\rangle$ via a $\Sigma_1$ ($\Pi_1$)
 formula of ${\mathcal L}_{\in}(P)$. A relation on $L_{\alpha}[A]$ is   $\Delta_1^{L_{\alpha}[A]}$ if it is both $\Sigma_1^{L_{\alpha}[A]}$
 and $\Pi_1^{L_{\alpha}[A]}$.

 For a set $X$, $|X|$ denotes the cardinality of $X$. For further unexplained notions and
 proofs see \cite[II. pp. 102--104]{devlin} or \cite{hajnal,levy}.

 \begin{proposition}\label{eigen}
 \begin{enumerate}
 \item $\alpha\leq\beta\Rightarrow L_{\alpha}[A]\subseteq L_{\beta}[A]$.
 \item $\alpha<\beta\Rightarrow L_{\alpha}[A]\in L_{\beta}[A]$.
\item $L_{\alpha}[A]$ is transitive.
\item $L[A]\,\cap\,\alpha=L_{\alpha}[A]\,\cap\,\alpha=\alpha$.
\item For $\alpha\geq\omega$, $|L_{\alpha}[A]|=|\alpha|$.
\item $L[A]\models \mathbf{ZF}$.
\item $\nu\mapsto L_{\nu}[A]$ is uniformly $\Delta_1\hochal$ for limits
$\lambda>\omega$.
\item $B=A\,\cap\,L[A]\,\Rightarrow\,L[A]=L[B]\;\wedge\;(V=L[B])\hocha$.
 \item There is a $\Sigma_1$ formula $\mathrm{wo}(x,y,z)$ such that
     $$\mathbf{KP}\vdash \mbox{``$\{\langle x,y\rangle\mid \mathrm{wo}(x,y,a)\}$ is a wellordering of $L[a]$''}$$
     and if $<_{L[A]}$ denotes the wellordering of $L[A]$ determined by $\mathrm{wo}$, then for any limit $\lambda>\omega$,
     $$<_{L[A]}\,\cap\,L[A]\times L[A]\mbox{ is }\Sigma_1\hochal.$$

     \item $L[A]$ is model of $\mathbf{AC}$.
     \item $\lambda>\omega\mbox{ limit }\wedge\,B=A\,\cap\,L_{\lambda}[A]\,\Rightarrow\,
     L_{\lambda}[A]=L_{\lambda}[B].$
\end{enumerate}
 \end{proposition}

 \section{Computability over $L[A]$}
 In this section we develop the recursion theory of
 partial $\Sigma_1^{ L[A]}$ functions, that is functions (not necessarily everywhere defined) whose graphs are $\Sigma_1^{ L[A]}$.
Below we shall write $L_{\alpha}[A]\models \varphi$
 rather than the more correct $\langle L_{\alpha}[A],\in,L_{\alpha}[A]\cap A\rangle \models \varphi$. Likewise, $\langle L[A],\in,L[A]\cap A\rangle \models \varphi$ will be shortened to $L[A]\models \varphi$.

 \begin{definition}{\em $\langle a,b\rangle$ denotes the ordered pair of
 two sets $a$ and $b$. If $c$ is an ordered pair $\langle a,b\rangle$ let $(c)_0=a$ and
 $(c)_1=b$. If $c$ is not an ordered pair let $(c)_0=(c)_1=0$.
 we also define ordered $n$-tuples via $\langle a_1\rangle:=a_1$ and
 $\langle a_1,\ldots,a_n,a_{n+1}\rangle:=\langle \langle a_1,\ldots,a_n\rangle,a_{n+1}\rangle$.

 It's standard procedure to assign to each formula $\psi$ of ${\mathcal L}_{\in}(P)$ a
 G\"odel number $\goed{\psi}$ such that $\goed{\psi}$ is a hereditarily definable
 set, for instance by using the pairing function $a,b\mapsto \langle a,b\rangle$.
 There is a formula $\mathrm{Sat}(v,w)$ of
 ${\mathcal L}_{\in}(P)$ such that for all $\Delta_0$ formulae $\theta(x_1,\ldots,x_n)$ of ${\mathcal L}_{\in}(P)$, not involving other free variables,
 the following holds 
 for any limit $\lambda>\omega$ and all
 $\vec a=a_1,\ldots,a_n\in L_{\lambda}[A]$:
 \begin{eqnarray}\label{rek1}  L_{\lambda}
 \models \theta(\vec a\,) &\mbox{ iff }&
  L_{\lambda}\models
 \mathrm{Sat}(\goed{\theta},\langle \vec a\,\rangle).\end{eqnarray}
 Moreover, $\mathrm{Sat}$ is uniformly $\Delta_1\hochal$ for
 limits $\lambda>\omega$ (see \cite[II]{devlin}).

 Now let $\lambda$ be a limit $>\omega$.
 For $e,a_1,\ldots,a_n\in L_{\lambda}[A]$ define
 \begin{eqnarray}\label{rek2}&&[e]\hochal_n(a_1,\ldots,a_n)\simeq b\end{eqnarray} if $e$ is an ordered pair $\langle\goed{\psi},c\rangle$ with $\psi$ being a $\Delta_0$ formula of ${\mathcal L}_{\in}(P)$,
 not involving free variables other than $x_1,\ldots,x_{n+2}$, such that
 \begin{eqnarray}\label{rek3}&& L_{\lambda}[A]\models
 \mathrm{Sat}(\goed{\psi},\langle a_1,\ldots,a_n,c,d\rangle)\end{eqnarray}
 and $(d)_0=b$, where $d$ is the $<_{L[A]}$-least ordered pair satisfying (\ref{rek3}).

 Likewise, $[e]\hocha_n(a_1,\ldots,a_n)\simeq b$ is defined by replacing $L_{\lambda}[A]$ by $L[A]$ in the foregoing definition.
 }\end{definition}

 \begin{lemma}\label{4.2} Let $\tau>\omega$ be a limit of limits, i.e.,
 $$\forall \xi<\tau\,\exists \lambda<\tau\,[\xi<\lambda\,\wedge\,\mbox{$\lambda$ limit}].$$
 \begin{itemize}
 \item[(i)] For $e\in L_{\tau}[A]$, the partial function $f$ on $L_{\tau}[A]$ given by
 \begin{eqnarray*}f(a_1,\ldots,a_n)=b &\mbox{ iff }& [e]^{L_{\tau}[A]}_n(a_1,\ldots,a_n)\simeq b\end{eqnarray*} is $\Sigma_1^{L_{\tau}[A]}$ (uniformly for all such $\tau$).
 \item[(ii)] For every $n$-ary partial $\Sigma_1^{L_{\tau}[A]}$ function $f$ there exists an index $e\in L_{\tau}[A]$ such that, for all $a_1,\ldots,a_n\in L_{\tau}[A]$,
 \begin{eqnarray*}f(a_1,\ldots,a_n)=b&\mbox{ iff }& [e]^{L_{\tau}[A]}_n(a_1,\ldots,a_n)\simeq b.\end{eqnarray*}

     \item[(iii)] (i) and (ii) hold with $L[A]$ in place of $L_{\tau}[A]$.
\item[(iv)]   $[e]^{L_{\tau}[A]}_n(a_1,\ldots,a_n)\simeq b$ implies $[e]\hocha_n(a_1,\ldots,a_n)\simeq b$ and $[e]\hochal_n(a_1,\ldots,a_n)\simeq b$ for all limits $\lambda>\tau$.
\item[(v)] If $[e]\hocha_n(a_1,\ldots,a_n)\simeq b$ then $[e]\hochal_n(a_1,\ldots,a_n)\simeq b$
for some limit $\lambda$.
\end{itemize}
\end{lemma}
\prf (i) First note that by Proposition \ref{eigen} the relation $<_{L[A]}$ restricted  to $L_{\lambda}[A]$ is $\Sigma_1\hochal$ for all limits $\lambda>\omega$. Thus the  $<_{L[A]}$-leastness of $d$ with respect to  (\ref{rek3})
can be expressed by
\begin{eqnarray*} &&\exists \lambda<\tau\,[\mbox{$\lambda$ limit $>\omega$}\,\wedge\,
a_1,\ldots,a_n,c,d \in L_{\lambda}[A]\\
&&\;\;\;\wedge\, L_{\lambda}[A]\models
\mbox{Sat}(\goed{\vartheta},\langle a_1,\ldots,a_n,c,d\rangle)\\
&&\;\;\; \forall u\in L_{\lambda}[A](u<_{L[A]}d\;\to\;
 L_{\lambda}[A]\models\mbox{Sat}(\goed{\neg \vartheta},\langle a_1,\ldots,a_n,c,u\rangle))],
\end{eqnarray*}
which is clearly $\Sigma_1^{L_{\tau}[A]}$.
\\[1ex]
(ii) Since $f$ is $\Sigma_1^{L_{\tau}[A]}$ there is a $\Sigma_1$ formula
$\exists x_{n+3}\vartheta_0(x_1,\ldots,x_{n+3})$ of ${\mathcal L}_{\in}(P)$
and a parameter $c\in \Sigma_1^{L_{\tau}[A]}$ (several  parameters can be coded as one) such that
 \begin{eqnarray*}f(a_1,\ldots,a_n)=b&\mbox{ iff }& L_{\tau}[A]\models \exists x_{n+3}\vartheta_0(a_1,\ldots,a_n,c,b,x_{n+3}).\end{eqnarray*}
 Now let $$\vartheta(x_1,\ldots,x_{n+2})\equiv \vartheta_0(x_1,\ldots,x_{n+1},(x_{n+2})_0,(x_{n+2})_1).$$
 Then
 \begin{eqnarray*}f(a_1,\ldots,a_n)=b&\mbox{ iff }&L_{\tau}[A]\models  \vartheta_0(a_1,\ldots,a_n,c,(d)_0,(d)_1]\,\mbox{ and }\,(d)_0=b,\end{eqnarray*}
 where $d$ is the $<_{L[A]}$-least $u$ such that
$ L_{\tau}[A]\models  \vartheta_0[a_1,\ldots,a_n,c,(u)_0,(u)_1]$.
Hence, with $e=\langle \goed{\vartheta},c\rangle$, we have $f(a_1,\ldots,a_n)=b$ iff $[e]_n^{L_{\tau}[A]}(a_1,\ldots,a_n)\simeq b$.

(iii) is proved in the same way as (i) and (ii).

(iv) follows since $\Sigma_1$ statements are upward persistent.

(v) follows since the statement is of $\Sigma_1$ form.
\qed

In several respects the recursion theory of partial $\Sigma_1\hocha$ functions and
partial $\Sigma_1^{L_{\tau}[A]}$ functions (for $\tau$ being a limit of limits)
shares similarities with ordinary recursion theory over $\omega$.
In particular, the analogues of the S-m-n theorem and the recursion theorem
hold.

 \section{Realizability over $L[A]$}
 $L[A]$ will be employed as a realizability universe. There is  a germane notion of realizability where realizers are indices of partial $\Sigma_1^{ L[A]}$ functions.
 \begin{definition}{\em For $d\in L[A]$ and set-theoretic sentences $\psi$ with
 parameters from $L[A]$ we define the realizability relation $d\Vdasha \psi$.

 Below we shall write $[e]\hocha(\vec a\,)\Vdasha \psi$
 rather than the more accurate $$\exists u\in L[A]([e]\hocha_n(\vec a\,)\simeq u\,\wedge\,u\Vdasha \psi]$$ where $\vec a=a_1,\ldots,a_n$.
 It will also assumed that all quantifiers range over $L[A]$.
 \begin{eqnarray*} e\Vdasha c\in d&\mbox{ iff }& c\in d\\
     e\Vdasha c= d&\mbox{ iff }& c= d\\
     e\Vdasha \varphi\wedge\psi &\mbox{ iff }& (e)_0\Vdasha \varphi\mbox{ and } (e)_1\Vdasha \psi\\
     e\Vdasha \varphi\vee\psi &\mbox{ iff }& [(e)_0=0\wedge (e)_1\Vdasha \varphi]\mbox{ or } [(e)_0=1\wedge (e)_1\Vdasha \psi]\\
      e\Vdasha \varphi\to\psi &\mbox{ iff }& \forall a\,[a\Vdasha \varphi\; \Rightarrow\; [e]\hocha( a)\Vdasha \psi]\\
       e\Vdasha \exists x\theta(x) &\mbox{ iff }& (e)_1\Vdasha \theta((e)_0)\\
       e\Vdasha \forall x\theta(x) &\mbox{ iff }& \forall a \, [e]\hocha(a)\Vdasha \theta(a).
       \end{eqnarray*}
      Occasionally we shall write $\Vdasha \psi $ for $\exists e\in L[A]\;e\Vdasha \psi$.
 }\end{definition}

 \begin{theorem}[Realizability Theorem]\label{realizability-thm}
 Let $\mathbb{R}\hocha$ be the set of real numbers in the sense of $L[A]$.
 If $\psi(x_1,\ldots,x_n)$ is a formula of set theory, with all free variables
 among the exhibited, and $\mathcal D$ is a proof of $ \psi(x_1,\ldots,x_n)$
 in $\mathbf T$, then one can effectively construct a hereditarily finite set
 $e_{\mathcal D}$ which only depends on $\mathcal D$ (and not on $A$) such that
 for all $a_1,\ldots,a_n\in L[A]$,
 \begin{eqnarray}\label{real1}
  && [e_{\mathcal D}]^{\mbox{\tiny $L[A]$}}(a_1,\ldots,a_n,\mathbb{R}\hocha)\Vdasha \psi(a_1,\ldots,a_n).\end{eqnarray}
  \end{theorem}
  \prf With little modification, the proof of Tharp's realizability theorem \cite{tharp} carries over to show this realizabilty theorem. This is similar to Tharp's realizability theorem \cite{tharp}. It can also be gleaned from the proofs
  of the realizability theorems \cite[Theorems 3.7-3.9]{schwach}, using considerable simplifications of
   the proofs brought about by the fact that there is uniform $\Sigma_1^{\mbox{\tiny $L[A]$}}$ selection function, i.e., there exists a hereditarily finite set
  $e_{ac}$ such that for all $A$ and nonempty sets $a$, $[e_{ac}]\hocha( a)\in a$.
   \qed

    \section{Designing $L[C]$}
    In order to show that $\CH\vee \neg \CH$ is not deducible in $\mathbf T$
    we intend to employ Theorem \ref{realizability-thm}. Aiming at a contradiction, we assume we have a derivation $\mathcal D$ of $\CH\vee \neg\CH$ in $\mathbf T$ and thus
    a hereditarily finite set $e_{\mathcal D}$ such that
    \begin{eqnarray}\label{Verm}&&[e_{\mathcal D}]({\mathbb R}\hocha)\Vdasha \CH\vee\neg\CH
    \end{eqnarray}
    holds for all sets $A$.
    To refute this, we intend  to carefully design a counterexample $C$.\footnote{Note that there are sets $A,A'$ and hereditarily finite sets $e,e'$ such that
    $[e]\hocha({\mathbb R}\hocha)\Vdasha \CH$ and
     $[e']\hochas({\mathbb R}\hochas)\Vdashas \neg\CH$, and hence sets
    $A$ such that $[e'']\hocha({\mathbb R}\hocha)\Vdasha CH\vee \neg CH$
    for some hereditarily finite $e''$.}
    We shall start from a set-theoretic universe $V_0$ such that
    $$V_0\models \mathbf{ZFC} +  2^{\aleph_0}=\aleph_2.$$
     $V_0$ can be obtained from any universe $V'$ such that $V'\models \mathbf{ZFC}+\mathrm{GCH}$ (e.g. $L$)
    by forcing with $\mathrm{Fn}(\kappa\times\omega,2)$, where the latter denotes the set of all finite functions with domain $\subset \kappa\times\omega$ and range $2$
    and
    $\kappa=(\aleph_2)^{V'}$, i.e., $\kappa$ is $\aleph_2$ in the sense of $V'$
    (see \cite[VII.5.14]{kunen}). Now let ${\mathbb R}^{\mbox{\tiny $V_0$}}$ be the reals in the sense of $V_0$. We would like to pick a set $C\in V_0$ such
    that ${\mathbb R}^{\mbox{\tiny $V_0$}}\in L[C]$. We cannot choose $C$ to be
    ${\mathbb R}^{\mbox{\tiny $V_0$}}$ since $L[{\mathbb R}^{\mbox{\tiny $V_0$}}]=L$ (cf. footnote 3)
    and therefore ${\mathbb R}^{\mbox{\tiny $V_0$}}\notin L[{\mathbb R}^{\mbox{\tiny $V_0$}}]$. But since $V_0$ satisfies $\AC$ there is an injective function $F$  in $V_0$ with domain ${\mathbb R}^{\mbox{\tiny $V_0$}}$ whose range is a set of ordinals.
    Identifying ${\mathbb R}^{\mbox{\tiny $V_0$}}$ with the set $\{g\in V_0\mid g:\mathbb N\to\mathbb N\}$, let
    \begin{eqnarray}\label{A} C &=& \{\omega^{F(g)+2}+\omega\cdot g(n)+n\mid g\in {\mathbb R}^{\mbox{\tiny $V_0$}}\}.
    \end{eqnarray} Then $C$ is a set of ordinals in $V_0$ and, owing to the uniqueness of the Cantor normal form,
     ${\mathbb R}^{\mbox{\tiny $V_0$}}$ is definable from $C$ in $L[C]$. The latter entails that  ${\mathbb R}^{\mbox{\tiny $V_0$}}\in L[C]$ and
     thus
     \begin{eqnarray}\label{reals} {\mathbb R}^{\mbox{\tiny $V_0$}}&=&
     {\mathbb R}^{\mbox{\tiny $L[C]$}}.\end{eqnarray}
     As a result, $L[C]\not\models \CH$ and therefore
     \begin{eqnarray}\label{w0} \mbox{for all $d\in L[C]$,
     $d\not \Vdasha \CH$.}\end{eqnarray} The assumption (\ref{Verm}) implies that
     there exists $b\in L[C]$ such that $L[C]\models [e_{\mathcal D}]({\mathbb R}\hochc)\simeq b$.
      Moreover,  (\ref{Verm}) and (\ref{w0}) entail that \begin{eqnarray}\label{w00} &&(b)_0=1.\end{eqnarray}
 We can now pick a sufficiently large limit ordinal $\rho$
      such that $C\in L_{\rho}[C]$, ${\mathbb R}^{\mbox{\tiny $V_0$}}\in L_{\rho}[C]$ and
      $b\in L_{\rho}[C]$. By Lemma \ref{4.2}(v) we can also arrange that
        \begin{eqnarray}\label{w1}&&L_{\rho}[C]\models [e_{\mathcal D}]({\mathbb R}\hochc)\simeq b.\end{eqnarray}
        Moreover, from Lemma \ref{4.2}(iv) and Proposition \ref{eigen}(11) it follows that for every
        set of ordinals $B$ with $B\cap\rho=\emptyset$ we have
         \begin{eqnarray}\label{w2}&&L[C\cup B]\models [e_{\mathcal D}]({\mathbb R}\hochc)\simeq b.\end{eqnarray}
         The next step consists in taking a forcing extension $V_1$ of $V_0$ which does not pick up new real
         numbers but satisfies $V_1\models 2^{\aleph_0}=\aleph_1$, i.e.,
         \begin{eqnarray}\label{w3} && {\mathbb R}^{\mbox{\tiny $V_0$}}={\mathbb R}^{\mbox{\tiny $V_1$}}\;\wedge\; (\aleph_1)^{\mbox{\tiny $V_0$}}=
  (\aleph_1)^{\mbox{\tiny $V_1$}}\;\wedge\;
          V_1\models\CH.\end{eqnarray}
         The latter can be arranged  by forcing with $${\mathbb P}:=(\mathrm{Fn}(\aleph_1,\aleph_2,\aleph_1))^{V_0}$$
         i.e., the set of functions $f\in V_0$ which are countable in $V_0$ with domain contained in $(\aleph_1)^{\mbox{\tiny $V_0$}}$ and range contained in $(\aleph_2)^{\mbox{\tiny $V_0$}}$.
        That (\ref{w3}) holds follows, e.g., from \cite[Ch.VII,6.13,6.14,6.15]{kunen}.

        Next we'd like to engineer a set $E\in V_1$ of ordinals all of whose members are greater than $\rho$ such that $L[ C\cup E]\models \CH$. Since $V_1$ is a model of the axiom of choice, there are functions $G$ and $H$ with domains $\{\alpha\mid \omega\leq \alpha<(\aleph_1)^{\mbox{\tiny $V_1$}}\}$ and $\{\beta\mid(\aleph_1)^{\mbox{\tiny $V_1$}}\leq \beta<(\aleph_2)^{\mbox{\tiny $V_1$}}\}$, respectively, such that for each $\alpha\in \dom(G)$, $G_{\alpha}:=G(\alpha)$ is a bijection between $\alpha$ and $\omega$, and for each $\beta\in\dom(H)$, $H_{\beta}:=H(\beta)$ is a bijection between $\beta$ and  $(\aleph_1)^{\mbox{\tiny $V_1$}}$. Let $\kappa$ and $\pi$ be fixed points of the function $\xi\mapsto \omega^{\xi}$ such that $\kappa<\pi$ and $\rho, (\aleph_1)^{\mbox{\tiny $V_1$}},(\aleph_2)^{\mbox{\tiny $V_1$}}<\kappa$.
        Now define
         \begin{eqnarray*}
           E_1 &:=& \{\kappa^{\alpha}\cdot(1+\xi)+G_{\alpha}(\xi)\mid \alpha\in\dom(G)\,\wedge\,\xi<\alpha\}\\
           E_2 &:=& \{\pi^{\beta}\cdot(1+\gamma)+H_{\beta}(\gamma)\mid \beta\in\dom(H)\,\wedge\,\gamma<\beta\}\\
           E&:=& E_1\cup E_2\end{eqnarray*}
           where of course $\kappa^{\alpha}$ and $\pi^{\beta}$ refer to the operation of ordinal exponentiation. Then $E_1\cap E_2=\emptyset$. Moreover, owing to the uniqueness of Cantor normal forms (e.g. \cite[Theorem 8.4.4]{takeuti}), for each $\alpha\in\dom(G)$, $F_{\alpha}$ is
           definable from $C\cup E$ in $L[C\cup E]$ (using the parameter $\kappa$), and likewise, for each $\beta\in \dom(H)$,
           $H_{\beta}$ is definable from $C\cup E$ in $L[C\cup E]$  (using the parameter $\pi$).
           To elaborate on this, suppose $\alpha\in\dom(G)$. Then for $\xi<\alpha$ search for the least ordinal $\delta$
           such that $\kappa^{\alpha}\cdot(1+\xi)+\delta\in E$. Necessarily, $\delta=G_{\alpha}(\xi)$.

           As a consequence of the above, we have
           \begin{eqnarray}\label{w4} && (\aleph_1)^{\mbox{\tiny $V_1$}} =(\aleph_1)^{\mbox{\tiny $L[C\cup E]$}}\;\wedge\;(\aleph_2)^{\mbox{\tiny $V_1$}} =(\aleph_2)^{\mbox{\tiny $L[C\cup E]$}}
           \;\wedge\;L[C\cup E]\models \CH.\end{eqnarray}
           To see the latter, suppose that $x\in L[C\cup E]$ and $x$ is an infinite set of reals.
           As $L[C\cup E]$ is a model of $\AC$, there is an ordinal $\eta$ and a bijection $\ell\in L[C\cup E]$ between $x$ and $\eta$. Since $L[C\cup E]\subseteq V_1$, $\eta<(\aleph_2)^{\mbox{\tiny $V_1$}}$ must obtain, and hence there is a bijection in $L[C\cup E]$ either between $\omega$ and $x$ or between $(\aleph_1)^{\mbox{\tiny $V_1$}} =(\aleph_1)^{\mbox{\tiny $L[C\cup E]$}}$ and $x$.
           From (\ref{reals}) and (\ref{w3}), we also conclude that
            \begin{eqnarray}\label{w5} && {\mathbb R}^{\mbox{\tiny $L[C]$}}={\mathbb R}^{\mbox{\tiny $L[C\cup E]$}}.\end{eqnarray}
            Utilizing the wellordering $<_{L[C\cup E]}$ and (\ref{w4}), there exists a $\Sigma_1^{\mbox{\tiny $L[C\cup E]$}}$ partial function $g$ that finds for each no-empty set of reals either
            a surjection of $\omega$ onto $x$ or a surjection of $x$ onto  ${\mathbb R}^{\mbox{\tiny $L[C]$}}$ since being such a mapping $f$ is a $\Delta_0$ property of $f$ in the parameters $x,\omega$ and
            ${\mathbb R}^{\mbox{\tiny $L[C]$}}$. Thus there is a realizer $d\in L[C\cup E]$ such that
            $d\Vdashce \CH$.
            From (\ref{w2}) and (\ref{Verm}) it then follows that $(b)_0=0$, contradicting (\ref{w00}).
            In sum, a contradiction has been inferred from (\ref{Verm}). On account of Theorem \ref{realizability-thm}, this implies that $\CH\vee\neg \CH$ is not provable in $\mathbf T$.

      \paragraph{Acknowledgement:} The work in this article was funded by a Leverhulme Trust Research Fellowship. The author is grateful to Sol Feferman for urging him to work on the problem.
      He is also grateful to Sy Friedman
      for answering questions about $L[A]$.



\begin{thebibliography}{TvD88}
 \bibitem{maralt} P. Aczel, M. Rathjen: {\em Notes on constructive set theory}, Technical Report 40,
Institut Mittag-Leffler (The Royal Swedish Academy of Sciences,
2001). {\tt http://www.mittag-leffler.se/preprints/0001/}, Preprint No. 40.

\bibitem{mar} P. Aczel, M. Rathjen: {\em Notes on constructive set theory},
Preprint (2010) 243 pages. http://www1.maths.leeds.ac.uk/~rathjen/book.pdf


 \bibitem{barwise} J. Barwise: {\em Admissible Sets and Structures}
(Springer, Berlin 1975).

\bibitem{beeson}
M.~Beeson: {\em Foundations of Constructive Mathematics}.
(Springer-Verlag, Berlin, Heidelberg, New York, Tokyo, 1985).

\bibitem{RMC} R.-M. Chen, M. Rathjen:  {\em Lifschitz Realizability for Intuitionistic Zermelo-Fraenkel  Set Theory.}
 Archive for Mathematical Logic 51 (2012) 789--818.

 \bibitem{devlin} K. Devlin: {\em Constructibility}. (Springer, Berlin, Heidelberg, New York, Tokyo, 1984).

 \bibitem{Diaconescu}    R. Diaconescu: {\em Axiom of choice and complementation}. Proc. Amer.
Math. Soc. 51 (1975) 176--178.

\bibitem{dummett} M. Dummett: {\em Frege: Philosophy of Mathematics} (Harvard University Press, London, 1991).


\bibitem{F1} S. Feferman: {\em On the strength of some semi-constructive theories}.
In: U.Berger, P. Schuster, M. Seisenberger (Eds.): {\em Logic, Construction, Computation} (Ontos Verlag, Frankfurt, 2012) 201--225.

\bibitem{F2} S. Feferman: {\em
Is the continuum hypothesis a definite mathematical
problem?}. Draft of paper for the lecture to the Philosophy Dept., Harvard University, Oct. 5, 2011 in the {\em Exploring the Frontiers of Incompleteness} project series,
Havard 2011--2012. 

\bibitem{F3} S. Feferman: {\em
Three Problems for Mathematics: Lecture 2: Is the Continuum Hypothesis a definite mathematical problem?}. Slides for inaugural Paul Bernays Lectures, ETH, Z\"urich, Sept. 12, 2012. 

\bibitem{F4} S. Feferman: {\em Why isn't the Continuum Problem on the Millennium ($\$$1,000,000) Prize list?}. Slides for CSLI Workshop on Logic, Rationality and Intelligent Interaction, Stanford, June 1, 2013.

    \bibitem{frege} G. Frege: {\em Die Grundlagen der Arithmetik}. (Verlag Wilhelm Koebner, Breslau, 1884).

    \bibitem{frieda} H.~Friedman: {\em Some applications of Kleene's method
for intuitionistic systems.}
In: A. Mathias and H. Rogers (eds.): {\em Cambridge Summer School in Mathematical Logic},
volume 337 of {\em Lectures Notes in Mathematics} (Springer, Berlin, 1973) 113--170.

\bibitem{hajnal} A. Hajnal:  {\em On a Consistency Theorem Connected with the Generalised Continuum Problem.}
Zeitschrift f\"ur Math. Logik 2 (1956) 131--136

\bibitem{kunen} K. Kunen: {\em Set theory} (North-Holland, Amsterdam, New York, Oxford, 1980).
\bibitem{levy} A. L\'evy:
1957 {\em Ind\'ependence Conditionnelle de V = L et d'Axiomes qui se Rattachent au Systeme
de M. G\"odel}. C.R. Acad. Sci. Paris 245 (1957) 1582--1583.

\bibitem{mc86} D.C.~McCarty: {\em Realizability and recursive set theory}, Annals of Pure and Applied
Logic 32 ,(1986) 153--183.

\bibitem{myhill73} J.~Myhill: {\em Some properties of Intuitionistic Zermelo-Fraenkel set theory.}
In: A. Mathias and H. Rogers (eds.): {\em Cambridge Summer School in Mathematical Logic},
volume 337 of {\em Lectures Notes in Mathematics} (Springer, Berlin, 1973) 206--231.

 \bibitem{pozsgay1} L. Pozsgay: {\em  Liberal intuitionism as a basis for set theory}, in Axiomatic Set Theory,
Proc. Symp. Pure Math. XIII, Part 1 (1971) 321-330.

\bibitem{pozsgay2} L. Pozsgay: {\em Semi-intuitionistic set theory}, Notre Dame J. of Formal Logic 13 (1972)
546-550.

\bibitem{R04} M.~Rathjen: {\em Realizability for constructive Zermelo-Fraenkel set
theory}. In: J. V\"a\"an\"anen, V. Stoltenberg-Hansen (eds.): {\em
Logic Colloquium 2003}. Lecture Notes in Logic 24 (A.K. Peters,
2006) 282--314.

\bibitem{ober} M.~Rathjen:  {\em The formulae-as-classes interpretation of constructive
set theory.} In:  H. Schwichtenberg, K. Spies (eds.): {\em Proof
Technology and Computation}
 (IOS Press, Amsterdam,2006) 279--322.

\bibitem{tklr} M.~Rathjen: {\em The disjunction and other properties for constructive
Zermelo-Fraenkel set theory.} Journal of Symbolic Logic
70 (2005) 1233--1254.

\bibitem{tklrac} M.~Rathjen: {\em Metamathematical
 Properties of Intuitionistic  Set Theories with Choice Principles.}
 In:  S. B. Cooper, B. L\"owe, A. Sorbi (eds.): {\em New
Computational Paradigms: Changing Conceptions of What is
Computable} (Springer, New York, 2008) 287--312.

\bibitem{schwach} M . Rathjen: {\em From the weak to the strong existence property,} Annals of Pure and Applied Logic
163 (2012) 1400–-1418.

\bibitem{takeuti} G. Takeuti, W.M. Wilson: {\em Introduction to axiomatic set theory}.
(Springer, New York, heidelberg, Berlin, 1971).

    \bibitem{tharp} L.~Tharp: {\em A quasi-intuitionistic set theory}. Journal of Symbolic Logic
36 (1971) 456--460.

\bibitem{thiele} E.J. Thiele: {\em \"Uber endlich axiomatisierbare Teilsysteme der Zermelo-Fraenkel'schen Mengenlehre},
Zeitschrift f\"ur Mathematische Logik und Grundlagen der Mathematik 14 (1968) 39-58.

\bibitem{wolf} R. S. Wolf: {\em Formally Intuitionistic Set Theories with Bounded Predicates
Decidable}, PhD Thesis (Stanford University, 1974).




\end{thebibliography}
\end{document}